\documentclass{llncs} %
\usepackage{epsfig}
\usepackage{graphicx}
\usepackage{epstopdf}
\usepackage{color,amsmath,amssymb}

\newcommand{\mE}{{\mathbb E}}
\begin{document}

\title{Optimal mass transport over bridges}
\author{Yongxin Chen\inst{1}, Tryphon Georgiou\inst{1} and Michele Pavon \inst{2}}
\institute{Department of Electrical and Computer Engineering,
University of Minnesota\\
200 Union street S.E. Minneapolis, Minnesota MN 55455, U.S.A.\\ \email{\{chen2468,tryphon\}@umn.edu} \and Dipartimento di Matematica,
Universit\`a di Padova\\ via Trieste 63, 35121 Padova, Italy\\\email{pavon@math.unipd.it}}

\maketitle

\begin{abstract}
We present an overview of our recent work on implementable solutions to the Schr\"{o}dinger bridge problem and their potential application to optimal transport and various generalizations.
\end{abstract}

\section{Introduction}\label{sec:introduction}
In a series of papers, Mikami, Thieullen and L\'eonard \cite{leo2,leo,Mik,mt,MT} have investigated the connections between the optimal mass transport problem (OMT) and the Schr\"{o}dinger bridge problem (SBP). The former may be shown to be the $\Gamma$-limit of a sequence of the latter, and thereby, SBP can be seen as a regularization of the OMT. Since OMT is well-known to be challenging from a computational viewpoint, this observation leads to the question of whether we can get approximate solutions to OMT via solving a sequence of SBPs. Both types of problem admit a control, fluid-dynamic formulation and it is in this setting that the connection between the two becomes apparent. There are, however, several difficulties in carrying out  this program:
\begin{itemize}
\item[i)] The solution of the SBP is usually not given in {\em implementable form};
\item[ii)]  SBP has been studied only for {\em non degenerate, constant diffusion coefficient} processes with control and noise entering through {\em identical channels} (this excludes most engineering applications);
\item[iii)] No SBP {\em steady-state} theory;
\item[iv)] No OMT problem with {\em nontrivial prior}.
\end{itemize}
In the past year, we have set out to partially remedy this situation \cite{GP}-\cite{CGP7}. We present here an overview of this work.

\section{Background}
\subsection{Optimal transport}
Consider the Monge-Kantorovich (OMT) problem \cite{Vil1,AGS,Vil2}
$$\inf_{\pi\in\Pi(\mu,\nu)}\int_{\bbbr^n\times\bbbr^n}c(x,y)d\pi(x,y)\enspace,
$$
where $\Pi(\mu,\nu)$ are ``couplings" of $\mu$ and $\nu$, and $c(x,y)=\frac12\|x-y\|^2$.

If $\mu$ does not give mass to sets of dimension $\le n-1$, by Brenier's theorem, there exists a unique optimal transport plan $\pi$ (Kantorovich) induced by a map $T$ (Monge), where $T=\nabla\varphi$, $\varphi$  is a {\em convex} function, $\pi=(I\times\nabla\varphi)\#\mu$, and $\nabla\varphi\#\mu=\nu$ where $\#$ indicates ``push-forward".
Assume from now on $\mu(dx)=\rho_0(x)dx$, $\nu(dy)=\rho_1(y)dy$. The static OMT above was given a dynamical formulation by Benamou-Brenier in \cite{BB}:
\begin{eqnarray}
\label{BB1}&&\inf_{(\rho,v)}\int_{\bbbr^n}\int_{0}^{1}\frac{1}{2}\|v(x,t)\|^2\rho(t,x)dtdx\enspace,\\&&\frac{\partial \rho}{\partial t}+\nabla\cdot(v\rho)=0\enspace,\label{BB2}\\&& \rho(0,x)=\rho_0(x), \quad \rho(1,y)=\rho_1(y)\enspace.\label{BB3}
\end{eqnarray}
\begin{proposition}Let $\rho^*(x,t)$ with $t\in[0,1]$ and $x\in \bbbr^n$, satisfy
\begin{equation}\label{optev}
\frac{\partial \rho^*}{\partial t}+\nabla\cdot(\nabla\psi\rho^*)=0, \quad \rho^*(x,0)=\rho_0(x)\enspace,
\end{equation}
where $\psi$ is a (viscosity) solution of the Hamilton-Jacobi equation
\begin{equation}\label{HJclass}
\frac{\partial \psi}{\partial t}+\frac{1}{2}\|\nabla\psi\|^2=0
\end{equation}
for some boundary condition $\psi(x,1)=\psi_1(x)$.
If $\rho^*(x,1)=\rho_1(x)$, then the pair $\left(\rho^*,v^*\right)$ with $v^*(x,t)=\nabla\psi(x,t)$ is a solution of (\ref{BB1})-(\ref{BB3}).
\end{proposition}

\subsection{Schr\"{o}dinger bridges}The ingredients of a Schr\"{o}dinger bridge problem are the following:
\begin{itemize}
\item a cloud of $N$ independent Brownian particles,
\item an initial and a final marginal density $\rho_0(x)dx$ and $\rho_1(y)dy$, resp.,
\item $\rho_0$ and $\rho_1$ are not compatible with the transition mechanism
$$\rho_1(y)\neq \int_{0}^{1}p(0,x,1,y)\rho_0(x)dx\enspace,
$$
where
$$p(s,y,t,x)=\left[2\pi(t-s)\right]
^{-\frac{n}{2}}\exp\left[-\frac{|x-y|^2} {2(t-s)}\right],\quad s<t\enspace.
$$
\end{itemize}
In view of the law of large numbers, particles have been transported in an {\em unlikely way} ($N$ being large). Then,
Schr\"{o}dinger in (1931) posed the following question: {\em Of the many unlikely ways in which this could have happened, which one is
the most likely?} F\"{o}llmer in 1988 observed that this is a problem of {\em large deviations of the empirical distribution} \cite{DZ} on path space  connected through Sanov's theorem to a {\em maximum entropy problem}.

Schr\"{o}dinger's solution ({\em bridge from $\rho_0$ to $\rho_1$ over Brownian motion}) has at each time a density $\rho$ that factors as $\rho(x,t)=\varphi(x,t)\hat{\varphi}(x,t)$, where $\varphi$ and $\hat{\varphi}$ solve the {\em Schr\"{o}dinger's system}
\begin{eqnarray}\label{1}
&&\varphi(x,t)=\int
p(t,x,1,y)\varphi(y,1)dy,\quad \varphi(x,0)\hat{\varphi}(x,0)=\rho_0(x)\enspace,\\&&\hat{\varphi}(x,t)=\int
p(0,y,t,x)\hat{\varphi}(y,0)dy,\quad \varphi(x,1)\hat{\varphi}(x,1)=\rho_1(x)\enspace.\label{2}
\end{eqnarray}
The new evolution has drift field $b(x,t)=\nabla\varphi(x,t)$. His result extends to the case when the ``prior" evolution is a general Markov diffusion process possibly with creation and killing \cite{W}. Existence and uniqueness for the Schr\"{o}dinger's system has been studied in particular by Beurling, Fortet, Jamison and F\"{o}llmer \cite{Beu,For,Jam,F2}, see \cite{W,leo} for a survey.

The {\em maximum entropy} formulation of the Schr\"{o}dinger bridge problem (SBP) with ``prior" $P$ is
$$\mbox{Minimize}\quad H(Q,P)=\mE_Q\left[\log\frac{dQ}{dP}\right] \quad \mbox{over} \quad \mathcal D(\rho_0,\rho_1)\enspace,
$$
where  $\mathcal D$ is the family of
distributions on $\Omega:=C([0,1],\bbbr^n)$ that are equivalent to stationary Wiener measure $W=\int W_x\,dx$.
It can be turned, thanks to {\em Girsanov's theorem},  into a stochastic control problem see \cite{DP,DPP,PW,FHS} with fluid dynamic counterpart. Here $P=W^{\epsilon}$, namely stationary Wiener measure with variance $\epsilon$, in which case the problem becomes
    \begin{eqnarray}\nonumber
        && \inf_{(\rho,v)} \int_{\bbbr^n}\int_{0}^{1}\frac{1}{2\epsilon}\|v(x,t)\|^2\rho(x,t)dtdx\enspace,
        \\
        && \frac{\partial \rho}{\partial t}+\nabla\cdot(v\rho)- \nonumber
        \frac{\epsilon}{2}\Delta \rho=0\enspace,
        \\
        && \rho(x,0)=\rho_0(x), \quad \rho(y,1)=\rho_1(y)\enspace.\nonumber
    \end{eqnarray}
This formulation should be compared to (\ref{BB1})-(\ref{BB3}).

\section{Gauss-Markov bridges}
Consider the problem in the case where the prior evolution and the marginals are {\em Gaussian}. In \cite{CGP1,CGP3}, the following two problems have been addressed:\\[.1in]
 \noindent
{\bf Problem 1:}
Find a control $u$, adapted to $X_t$ and minimizing
\[
J(u):=\mE\left\{\int_0^1u(t)\cdot u(t) \,dt\right\}\enspace,
\]
among those which achieve the transfer
\begin{eqnarray}\nonumber
dX_t=A(t)X_tdt+B(t)u(t)dt+B_1(t)dW_t\enspace,\\ X_0\sim\mathcal N(0,\Sigma_0),\quad X_1\sim\mathcal N(0,\Sigma_1)\enspace.\nonumber
\end{eqnarray}

If the pair $(A,B)$ is controllable (for constant $A$ and $B$, this amounts to the matrix $\left(B,AB,...,A^{n-1}B\right)$ having full row rank), Problem 1 turns out to be always {\em feasible} (this result is highly nontrivial as the control may be ``handicapped" with respect to the effects of the noise).\\[.1in]
\noindent
{\bf Problem 2:}
Find $u=-Kx$ which minimizes $J_{\mathrm{power}}(u):=\mE \{u\cdot u\}$ and such that
$$
dX_t=(A-BK)X_tdt+B_1dW_t
$$
has
$$\rho(x)=(2\pi)^{-n/2}\det (\Sigma)^{-1/2}\exp\left(-\frac{1}{2}x'\Sigma^{-1}x\right)
$$
as {\em invariant probability density}.

Problem 2 may not have a solution (not all values for $\Sigma$ can be maintained by state feedback).

\noindent
{\em Sufficient conditions} for optimality have been provided in \cite{CGP1,CGP3} in terms of:
\begin{itemize}
\item{}  a system of two {\em matrix Riccati equations} ({\em Lyapunov equations} if $B=B_1$) in the finite horizon case. The Riccati equations are nonlinearly coupled through the boundary conditions. In the case where $B\neq B_1$, which falls outside the classical maximum entropy problem but represents a relaxed version of the classical steering problem, the two equations are also {\em dynamically} coupled.
\item{}in terms of {\em algebraic conditions} for the stationary case.
\end{itemize}
\noindent
Optimal controls may be computed via semidefinite programming in both cases.

\section{Cooling for stochastic oscillators}
Cooling for micro and macro-mechanical systems consists in implementing via feedback a frictional force to steer the state of a thermodynamical system to a non equilibrium steady state with {\em effective temperature} that is lower than that of the heat bath. Important applications of such {\em Brownian motors} \cite{R} are found in molecular dynamics, Atomic Force Microscopy and gravitational wave detectors \cite{DE,LMC,Vin}, to name a few.

The basic model is provided by a {\em controlled stochastic oscillator} deriving from the Nyquist-Johnson model of RLC electrical network with noisy resistor (1928) and the Ornstein-Uhlenbeck model of physical Brownian motion (1930):
\begin{eqnarray}\label{OU1}
dx(t)&=&v(t)\,dt\enspace,\\\label{OU2}
dv(t)&=&-\beta v(t)\,dt-\frac{1}{m}\nabla V(x(t))dt+ u(x(t), v(t),t)+\sigma dW(t)\enspace,\\\label{OU3}\sigma^2&=&\frac{2k\beta T}{m},\quad \mbox{\em Einstein's fluctuation-dissipation relation}.
\end{eqnarray}
Here $u(x,v,t)$ is a feedback control law and $V$ is such that the initial value problem is well-posed on bounded time intervals. For $u\equiv 0$,
$$\rho(x,v,t)\rightarrow\rho_{MB}(x,v)=Z^{-1}\exp\left[-\frac{H(x,v)}{kT}\right],\; H(x,v)=\frac12 m v\cdot v +V(x)\enspace.
$$
Let $\bar{\rho}(x,v)=C\exp\left[-\frac{H(x,v)}{kT_{{\rm eff}}}\right]$ and let $T_{{\rm eff}}<T$ be a desired {\em steady state} effective temperature. In \cite{CGP4}, we have studied the following two problems:
\begin{itemize}
\item Efficient asymptotic steering of the system to $\bar{\rho}$;
\item Efficient steering of the system from the initial condition $\rho_0$ to $\bar{\rho}$ at a finite time $t=1$.
\end{itemize}
In both cases, we get a solution for a general system of nonlinear stochastic oscillators, where we allow for both potential and dissipative interactions between the particles, by extending the theory of the Schr\"{o}dinger bridges accordingly.

Consider the case of a scalar oscillator in a quadratic potential with Gaussian marginals. For a suitable choice of constants, the model is
\begin{eqnarray*}
dx(t)&=&\phantom{-}v(t)dt\enspace,\\
dv(t)&=& -v(t)dt -x(t)dt+u(t)dt +dW(t)\enspace.
\end{eqnarray*}
Using the results in \cite{CGP1,CGP3}, through velocity feedback control the system is first efficiently steered to the desired state state $\bar{\rho}$ at time $t=1$ and then maintained efficiently in $\bar{\rho}$. This is illustrated by Figure $1$ that depicts some sample paths and a transparent tube outlining the ``$3\sigma$ region'' of the one-time densities.

\begin{figure}[htb]
\centering
\includegraphics[totalheight=2in]{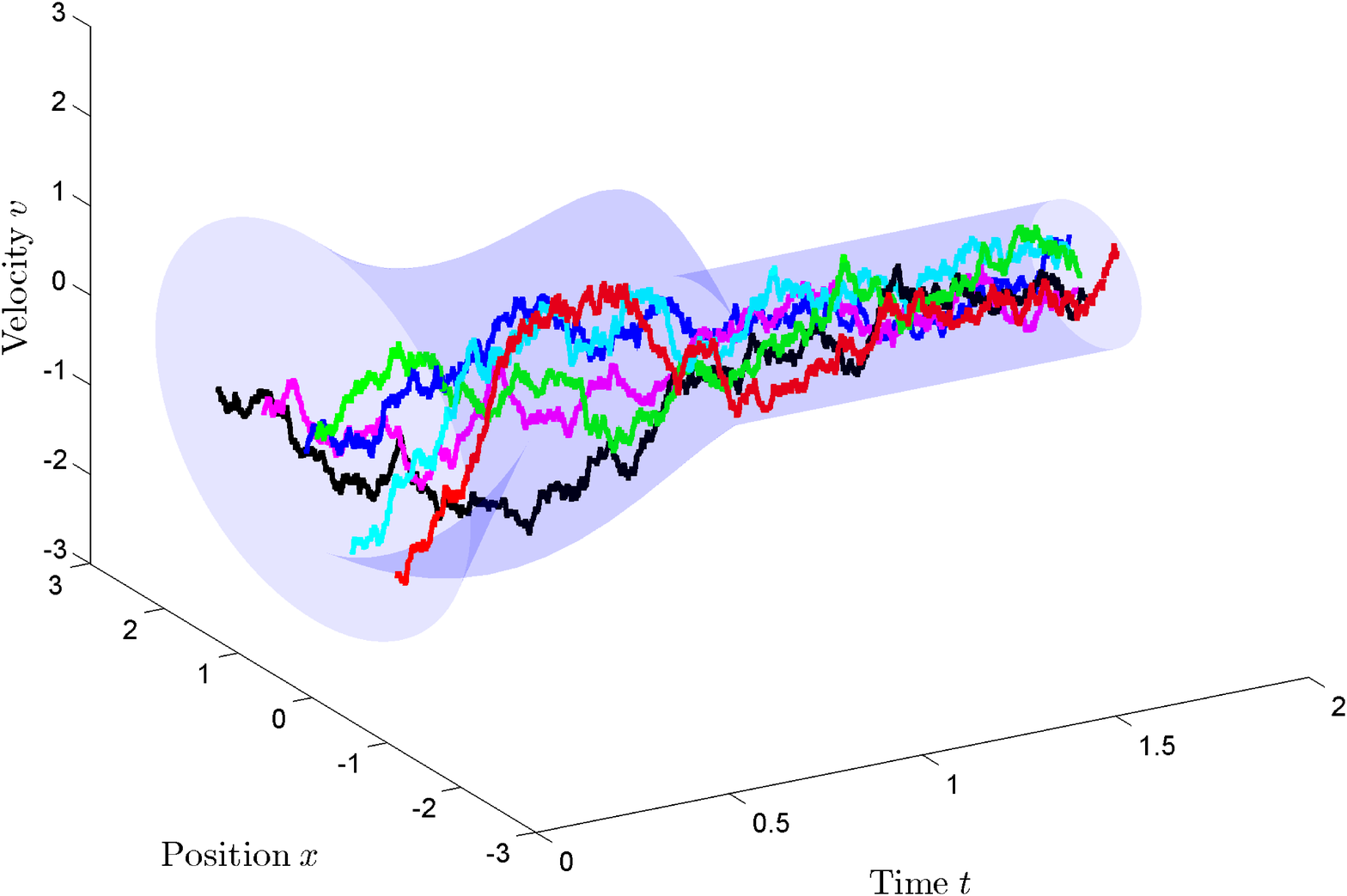}
\caption[ ]{Inertial particles: trajectories in phase space.}
\end{figure}

\section{OMT with prior}
In \cite{CGP5,CGP6}, we have formulated and studied a generalization of optimal transport problem that includes prior dynamics. It is the natural candidate for the zero-noise limit of SBP where the prior is a general Markovian evolution and not just stationary Wiener measure. In particular, in \cite{CGP6} we have studied the case where there are fewer control than state variables and Gaussian marginals and derived the corresponding limiting transport problem. The latter can be put in the form of a classical OMT with cost deriving from a Lagrangian action, where, however, the Lagrangian is not strictly convex with respect to the $\dot{x}$ variable. Convergence of solutions is proven directly. Simulations confirm that in the zero-noise limit the ``entropic interpolation" provided by the (generalized) Schr\"{o}dinger bridge converges to the ``displacement interpolation" of the limiting OMT problem.

In conclusion, in \cite{CGP1,CGP2,CGP4,CGP5,CGP6}, we have worked out a number of cases where an implementable form of the solution of a (possibly generalized) Schr\"{o}dinger bridge problem can be obtained. We have also explored to some extent the connection between zero-noise limits of SBP and suitable reformulations of OMT problems. These cases include degenerate, hypoelliptic diffusions like the Ornstein-Uhlenbeck model (\ref{OU1})-(\ref{OU3}). The case of differing noise and control channels which does not have a classical SBP counterpart has also been studied. Finally, in \cite{CGP3}, we have extended the fluid-dynamic SBP theory to the case of anisotropic diffusions with killing, a situation where again no probabilistic counterpart is available in general. The new evolution is obtained by solving a suitable generalization of the Schr\"{o}dinger bridge system. How can we solve this generalized Schr\"{o}dinger system as well as those corresponding to problems not covered in  \cite{CGP1,CGP2,CGP4,CGP5,CGP6}? An alternative powerful tool is given by iterative schemes which contract Birkhoff's version of Hilbert's metric. This is discussed in the next section.
\section{Positive contraction mappings for Schr\"{o}dinger systems}

Let $\mathcal S$ be a real Banach space and $\mathcal K$ a closed solid cone in $\mathcal S$. That is, $\mathcal K$ is closed with nonempty interior and is such that $\mathcal K+\mathcal K\subseteq \mathcal K$, $\mathcal K\cap -\mathcal K=\{0\}$ as well as $\lambda \mathcal K\subseteq \mathcal K$ for all $\lambda\geq 0$. Define $x\preceq y \Leftrightarrow y-x\in\mathcal K$, and for $x,y\in\mathcal K\backslash \{0\}$, $M(x,y):=\inf\, \{\lambda\,\mid x\preceq \lambda y\}$ and $m(x,y):=\sup \{\lambda \mid \lambda y\preceq x \}$.
The {\em Hilbert metric} is the projective metric defined on $\mathcal K\backslash\{0\}$ by
\[
d_H(x,y):=\log\left(\frac{M(x,y)}{m(x,y)}\right).
\]
A map $\mathcal E$ from $\mathcal S$ to $\mathcal S$ is said to be {\em positive} provided it takes the interior of $\mathcal K$ into itself. For such a map define its {\em projective diameter}
\begin{eqnarray*}
\Delta(\mathcal E):=\sup\{d_H(\mathcal E(x),\mathcal E(y))\mid x,y\in \mathcal K\backslash\{0\}\}
\end{eqnarray*}
and the {\em contraction ratio}
\begin{eqnarray*}
\|\mathcal E\|_H:=\inf\{\lambda \mid d_H(\mathcal E(x),\mathcal E(y))\leq \lambda d_H(x,y),\mbox{ for all }x,y\in\mathcal K\backslash\{0\}\}.
\end{eqnarray*}

\begin{theorem} (Garret Birkhoff 1957, P. Bushell 1973)
{Let $\mathcal E$ be a positive map. If $\mathcal E$ is monotone and homogeneous of degree $m$ ($\mathcal E(\lambda x)=\lambda^m \mathcal E(x)$),
then it holds that}
\[
\|\mathcal E\|_H\leq m.
\]
If $\mathcal E$ is also linear, the (possibly stronger) bound also holds
\[
\|\mathcal E\|_H=\tanh(\frac{1}{4}\Delta(\mathcal E)).
\]
\end{theorem}

Consider now a Markov chain with $T$-step transition probabilities $\pi_{x_0,x_T}$ (prior) and consider two marginal distributions ${\mathbf p}_0$ and ${\mathbf p}_T$, where $x_0,x_T$ are indices corresponding to initial and final states.
An adaptation of Schr\"odinger's question to this setting leads to the following Schr\"{o}dinger system:
\begin{eqnarray}\nonumber
\varphi(0,x_0)&=&\sum_{x_T}\pi_{x_0,x_T}\varphi(T,x_T)=\mathcal E\left(\varphi(T,x_T)\right),\quad \varphi(0,x_0)\hat\varphi(0,x_0)={\mathbf p}_0(x_0),\\
\hat\varphi(T,x_T)&=&\sum_{x_0}\pi_{x_0,x_T}\hat\varphi(0,x_0)=\mathcal E^\dagger\left(\hat\varphi(0,x_0)\right),\quad \varphi(T,x_T)\hat\varphi(T,x_T)={\mathbf p}_T(x_T).\nonumber\end{eqnarray}
It turns out that the composition of the four maps
\begin{eqnarray}\nonumber
\hat\varphi(0,x_0) \longrightarrow \hat\varphi(T,x_T):=\mathcal E^\dagger(\hat\varphi(0,x_0))\longrightarrow\varphi(T,x_T):=\frac{{\mathbf p}_T(x_N)}{\hat\varphi(T,x_T)}\\\longrightarrow \varphi(0,x_0):=\mathcal E\left(\varphi(T,x_T)\right)\longrightarrow \left(\hat\varphi(0,x_0)\right)_{\rm next}:= \frac{{\mathbf p}_0(x_0)}{\varphi(0,x_0)}\nonumber
\end{eqnarray}
where division of vectors is performed componentwise, is {\em contractive in the Hilbert metric}. Indeed, the linear maps are non-expansive with $\mathcal E$ strictly contractive, whereas componentwise divisions are isometries (and contractive when the marginals have zero entries). In \cite{GP}, we have obtained similar results for Kraus maps of statistical quantum mechanics with pure states or uniform marginals. The case of diffusion processes is studied in \cite{CGP7}.  Applications include interpolation of 2D images to construct a 3D model (MRI).

 \end{document}